\input diagrams.tex


\centerline{\bf A homotopical algebra of graphs related to zeta series.}
\smallskip
\centerline{(Preliminary draft, February, 2008)}
\vskip .25 true cm
\centerline{Terrence Bisson \quad\& \quad Aristide Tsemo}
\centerline{bisson@canisius.edu \quad \  \quad tsemo58@yahoo.ca}
\vskip .25 true cm

\bigskip\noindent
{\bf Abstract:} The purpose of this paper is to
develop a homotopical algebra for graphs,
relevant to zeta series
and spectra of finite graphs. 
More precisely, we define a Quillen model structure
in a category of graphs
(directed and possibly infinite, 
with loops and multiple arcs allowed).
The weak equivalences for this model structure are the
Acyclics (graph morphisms which preserve cycles).
The cofibrations and fibrations for the model
are determined from the class of Whiskerings 
(graph morphisms produced by grafting trees).
Our model structure seems to fit well with
the importance of acyclic directed graphs in many applications.

In addition to the weak factorization systems 
which form this model structure,
we also describe two Freyd-Kelly factorization systems based on 
Folding, Injecting, and Covering graph morphisms.

\bigskip\noindent
{\bf 0. Introduction.}

In this paper we develop a notion of homotopy
within graphs, and demonstrate
its relevance to the study of zeta series
and spectrum of a finite graph.  
We will work throughout with a particular category of graphs,
described in Section 1 below.
Our graphs will be directed and possibly infinite, 
with loops and multiple arcs allowed.

\medskip
Let us explain what we mean by homotopy here.

We are not concerned with the geometric realization
of graphs as one-dimensional topological spaces.
Since one-dimensional CW complexes are homotopic 
to disjoint unions of joins of circles,
the usual invariants from algebraic topology 
cannot see much of the structure of a graph in this way.
In any case, {\it directed} graphs are definitely not just part of topology
(they are perhaps more related to new areas of {\it directed topology},
as in Fajstrup and Rosick\'y [2007]).

Homotopy originally referred to topological deformation of structure.
But Quillen's remarkable notes on homotopical algebra [1967] gave
abstract axioms for working with concepts of homotopy in rather general categories.
When these axioms are satisfied in a category, we say that we
have given a ``model structure'' there.
Quillen's axioms have led to new insights and developments
in settings such as
chain complexes and homological algebra,
simplicial sets, topos theory, and small categories 
(including monoids, groups, groupoids, and posets).
References include Thomason [1980], Joyal and Tierney [1991], 
Dwyer and Spalinski [1995], Cisinski [2002],  and many others.
Also, recent proofs of the Bloch-Kato and Milnor conjectures are based upon
development of a homotopical algebra for schemes;
see Voevodsky and Morel [1999]. 

\medskip

A central part of giving a model structure in a category is 
the specification of which morphisms
in the category are to be called ``weak equivalences''.
In most applications, the weak equivalences are defined to be those morphisms which
preserve some interesting invariant, such as  
homotopy type for topology, homology for chain complexes,
geometric realization for simplicial sets, and nerve or topos of presheafs 
for small categories.

In our model structure 
we use the cycle structure of directed graphs to determine our weak equivalences.
More precisely, we take as our weak equivalences the ``acyclic'' graph morphisms,
which neither create nor destroy cycles.
We hope that our model structure fits well with
the  role that acyclic directed graphs play in
in applications such as computer algorithms, analysis of the internet, 
random walks and markov chains, and representations of quivers.

In section 1 we set up our category Gph of graphs.  
In section 2 we give background on weak factorization systems in general,
and establish an example with classes of graph morphisms which we call Whiskerings and Surjectings.
In section 3, which is somewhat of a digression,
we discuss Freyd-Kelly factorization systems in general,
and consider two examples in Gph, 
involving classes of graph morphisms which we call Injectings, Foldings, and Coverings. 
In section 4 we give axioms for model structures in general, and 
define our model structure on Gph.
In section 5 we associate to each finite directed graph $X$ a zeta series $Z_X(u)$.
We show that if $f:X\to Y$ is an acyclic graph morphism, then $Z_X(u)=Z_Y(u)$
and the eigenvalues of the adjacency matrices of $X$ and $Y$ agree ``up to zero eigenvalues''.
The paper ends with two appendices.

\medskip
Acknowledgements: the first author would like to thank
Andr\'e Joyal and Ji\v r\'\i \ Rosick\'y for encouragement,
York University for hospitality during fall 2005, and 
his students in the ``Geometry and Physics on Graphs'' 
Research Experience for Undergraduates at Canisius College.

This paper uses Paul Taylor's excellent package diagrams.tex.

\bigskip\noindent
{\bf 1. A category of graphs.}

Let us establish precisely the objects and morphisms for our category Gph of graphs.
For our purposes,
a {\it graph} is a data structure $X=(X_0,X_1,s,t)$ with a set $X_0$ of {\it nodes},
a set $X_1$ of {\it arcs}, and a pair of functions $s,t:X_1\to X_0$ which
specify the {\it source} and {\it target} nodes of each arc.
We may say that $a\in X_1$ is an arc {\it from node $s(a)$ to node $t(a)$}.
A {\it graph morphism} $f:X\to Y$ is a pair of functions 
$f_1:X_1\to Y_1$ and $f_0:X_0\to Y_0$ such that
$s\circ f_1=f_0\circ s$ and $t\circ f_1=f_0\circ t$.

Let ${\bf N}$ denote the graph with one node and no arcs.
Let ${\bf A}$ denote the graph with one arc and two nodes (its source and target). 
Then the set of nodes of a graph $X$ 
can be identified with the set of graph morphisms from ${\bf N}$ to $X$,
and the set of arcs of $X$ can be identified with 
the set of graph morphisms from ${\bf A}$ to $X$.

There is a stimulating discussion in Lawvere [1989] of
Gph as the category of presheafs on the small category with objects $0$ and $1$ 
and two non-identity morphisms from $0$ to $1$.
It follows that ${\rm Gph}$ is a topos, and thus a category with many nice geometric and
algebraic and logical properties; see Mac Lane and Moerdijk [1994], for instance.
We just call attention to a few aspects here.

\medskip
The category Gph has all products, and all coproducts (sums);
it also has pull-backs (fiber products) and pushouts.
Also, products distribute over coproducts, etc.
As in any presheaf category,
these categorical constructions are performed ``elementwise'',
where a graph has two types of elements, the nodes and the arcs. 
For instance, the empty product (terminal object) $1$ is the graph
with one node and one arc (which is a loop); and
the empty coproduct (initial object) $0$ is the graph with no nodes and no arcs.   

The category Gph has geometric aspects; for instance, 
it is a cartesian closed category like a good category of ``spaces''.  
The category Gph also has logical aspects; for instance,
there is a graph $\Omega$ which acts as generalized truth-values
for graphs, in that graph morphisms $X\to\Omega$ classify sub-graphs of $X$
(see Session 32 in Lawvere and Schanuel [1997]).

\bigskip\noindent
{\bf 2. Two classes of graph morphisms, and a weak factorization system.}

The {\it path graph} ${\bf P}_n$ has nodes $\{i:0\leq i\leq n\}$
and arcs $\{ (i-1,i):1\leq i\leq n \}$, with $s((i-1,i))=i-1$ and $t((i-1,i))=i$.
Note that ${\bf P}_0={\bf N}$ and ${\bf P}_1={\bf A}$.
A {\it path of length $n$} in a graph $X$ is just a
graph morphism $\alpha:{\bf P}_n\to X$; we define $s(\alpha)=\alpha(0)$ and 
$t(\alpha)=\alpha(n)$.  If $a$ is an arc in $X$ such that $t(\alpha)=s(a)$,
 then we define $\alpha a:{\bf P}_{n+1}\to X$, the {\it concatenation of $\alpha$ and $a$}.

Let us introduce some useful shorthand for arcs in a graph $X$. 
For any node $x$ in a graph $X$, let $X(x,*)$ denote the set of those arcs 
in $X$ which have source $x$, and let $X(*,x)$ denote the set of arcs with target $x$.
Note that a graph morphism $f:X\to Y$ induces a function $f:X(x,*)\to Y(f(x),*)$, etc.

Here is a our first class of graph morphisms.

\medskip\noindent
{\bf Definition:} A graph morphism $f:X\to Y$ is {\it Surjecting} 
if $f:X(x,*)\to Y(f(x),*)$ 
is a surjective function for all $x\in X_0$.

\medskip 
A {\it discrete graph} is one with no arcs.
We say that a node $x$ is a {\it root} of the graph $X$ if $X(*,x)$ is empty.
Let $R(X)$ denote the set of roots in $X$, viewed as a discrete subgraph of $X$.
A {\it rooted tree} is a graph $T$ with one root $r$ such that,
for each each node $x$ in $T$, there is a unique (directed) path in $T$ from $r$ to $x$.
For example, the path graph ${\bf P}_n$ is a rooted tree; and
so is the infinite path ${\bf P}_\infty$ whose nodes are the set of non-negative integers and
whose arcs are the set of ordered-pairs $(i-1,i)$ of non-negative integers, 
with $s((i-1,i))=i-1$ and $t((i-1,i))=i$.

For any node $x$ in a graph $X$,
we can define a rooted tree $T_xX$, the {\it tree of paths in $X$ leaving $x$}.
The nodes in $T_xX$ are the finite paths in $X$ with source $x$
(note that $x$ is considered as a path of length $0$ in $X$); 
the arcs in $T_xX$ are the triples
$(\alpha,a,\alpha a)$ where $\alpha a$ is the concatenation of path $\alpha$
and arc $a$ in $X$; and $s(\alpha)=s(\alpha,a,\alpha a)$
and $t(\alpha,a,\alpha a)=t(a)$.
There are natural graph morphisms $T_xX\to X$ given by $\alpha\mapsto t(\alpha)$ and
$(\alpha,a,\alpha a)\mapsto a$. 

A {\it rooted forest} $F$ is a coproduct of rooted trees.   
If $F$ is a forest with roots $R(F)$, then we may form a new graph
$X_F$ as the pushout of the graph morphisms $r:R(F)\to F$ and $f:R(F)\to X$;
we say that $X_F$ is {\it formed by attaching the forest $F$ to $X$ along $R$}. 
For instance, if $r$ is a root in a tree $T$ and $x$ is any node in a graph $X$, 
then $X_T$ is formed by attaching the tree $T$ to $X$ at the node $x$
(along the graph morphisms $r:{\bf N}\to T$ and $x:{\bf N}\to X$).

\medskip\noindent
{\bf Definition:} A graph morphism $f:X\to Y$ is a {\it Whiskering} if 
$Y$ is formed by attaching some rooted forest to $X$.

For example, the graph morphism ${\bf s}:{\bf N}\to {\bf A}$ is a Whiskering,
where ${\bf s}$ exhibits the node graph ${\bf N}$ as the source subgraph 
of the arc graph ${\bf A}$.
Also, every isomorphism is a Whiskering, and 
$R(F)\to F$ is a Whiskering if $F$ is a rooted forest.

Our goal in this section is to demonstrate some remarkable factorization properties 
of Surjectings and Whiskerings. Here is the conceptual background. 

\medskip\noindent
{\bf Definition:} Let $\ell:X\to Y$ and $r:A\to B$ be morphisms in a category ${\cal S}$.
We say that $\ell\dagger r$ when, for all $f$ and $g$,
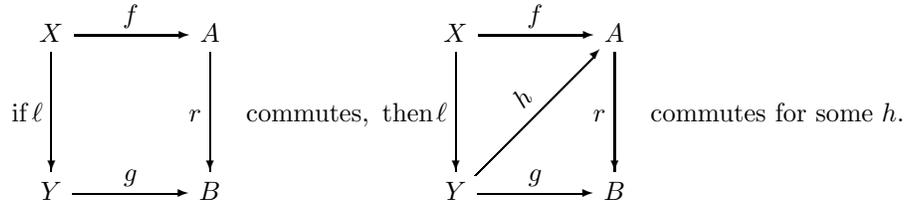
$${\rm if} \diagram
          X & \rTo^{f} & A  \cr
   \dTo^{\ell} &           &  \dTo^{r}   \cr
          Y  & \rTo^g   &B
                                    \enddiagram \quad{\rm commutes,\ then}
           \diagram
          X & \rTo^{f} & A  \cr
   \dTo^{\ell} &  \NE^{h} &  \dTo^{r}   \cr
          Y & \rTo^g   &B
                                    \enddiagram\quad {\rm commutes\ for\ some\ } h.$$ 
We say that $h$ is a filler for the commutative diagram.  We may also
say that $h$ {\it lifts} $g$ along $r$, or that $h$ {\it drops} (``extends'') $f$ along $\ell$.
Given two classes ${\cal L}$ and ${\cal R}$ of morphisms, we say 
${\cal L}\ \dagger\ {\cal R}$ when we have $\ell\dagger r$ for every $\ell\in{\cal L}$
and every $r\in {\cal R}$.
Given a class ${\cal F}$ of morphisms we may define
$${\cal F}^\dagger=\{r:f\dagger r, \ \forall f\in{\cal F}\}\quad{\rm and}\quad
{}^\dagger{\cal F}=\{\ell: \ell\dagger f, \ \forall f\in{\cal F}\}.$$

\medskip\noindent
{\bf Definition:} A {\it weak factorization system} in ${\cal S}$ 
is given by two classes
${\cal L}$ and ${\cal R}$ such that ${\cal L}^\dagger={\cal R}$ and
${\cal L}={}^\dagger{\cal R}$ and such that, for any morphism $c$ in ${\cal S}$,
there exist $\ell\in{\cal L}$ and $r\in{\cal R}$ with $c=r\circ\ell$.

The notion of weak factorization system has become a part of homotopical algebra;
see Section 4 here.  The appendix on Galois connections also provides some context.

\bigskip

The following three propositions combine to show that
Surjectings and Whiskerings give a weak factorization system in Gph.
Our Proposition 1 was inspired by an argument in Enochs and Herzog [1999].

\medskip\noindent
{\bf Proposition 1.} Any graph morphism $f:X\to Y$ factors as a Whiskering followed
by a Surjecting.

\medskip\noindent
{\bf Proof:} Recall that for each node $y$ in $Y$ we have
the tree $T_yY$ of paths in $Y$ leaving $y$.
From $f$ we construct the rooted forest $F=\sum_{x\in X_0} T_{f(x)}Y$, 
with roots $R=X_0$ considered as discrete subgraph of $X$.
The pushout of $R\to F$ along the subgraph inclusion $R\to X$ defines a
Whiskering $w:X\to X_F$.  We have $g:F\to Y$ as coproduct of the morphisms
$T_{f(x)}Y\to Y$, and since $f:X\to Y$ and $g:F\to Y$
agree on $R$, they determine a unique graph morphism $p:X_F\to Y$.
Note that $f=p\circ w$.

Let us show that $p$ is a Surjecting.  For any node $z$ in $X_F$ we must show
that $p:X_F(z,*)\to Y(p(z),*)$ is surjective.
But $z$ is either a node $x$ or a path in $Y$ with source $f(x)$,
for some $x\in X_0$.  

In the first case, we have $p:X_F(x,*)\to Y(f(x),*)$ with $X_F(x,*)=X(x,*)\cup T_yY(f(x),*)$,
and $T_yY(f(x),*)\to Y(f(x),*)$ is a bijection.

In the second case, we have $p:X_F(\alpha,*)\to Y(p(\alpha),*)$ with
$X_F(\alpha,*)=T_{f(x)}Y(\alpha,*)=Y(t(\alpha),*)=Y(p(\alpha),*)$.

In either case, $p:X_F(z,*)\to Y(p(z),*)$ is surjective.
QED

\bigskip\noindent
{\bf Proposition 2.} ${\rm Whiskering} \dagger {\rm Surjecting}$.

\medskip\noindent
{\bf Proof:} Let $f:Z\to Y$ be Surjecting.
First we show lifting of rooted trees.
More precisely, if $T$ is a rooted tree with root $x$ 
and we have the following commutative diagram
$$\diagram
    {\bf N} & \rTo^{z}  &Z     \cr
     \dTo^x &                &\dTo^{f}   \cr
          T & \rTo^{g}       & Y     
                                                    \enddiagram$$
then there is a filler $h:T\to Z$.  This follows by induction on the length of path
from root to nodes of $T$, as follows.
Suppose that we have extended $h$ to paths of length $n$ and let $\alpha a$ be a path
of length $n+1$ in $T$.  Let $x'=t(\alpha)$. Then $f(h(x'))=g(x')$ and $g(a)\in Y(g(x'),*)$
and $f:Z(h(x'),*)\to Y(g(x'),*)$ is a surjective function, so
there exists an arc $a'\in Z(h(x'),*)$ so that $f(a')=g(a)$.  We extend $h$ to $\alpha a$ by $h(a)=a'$.

More generally, consider any commutative diagram
$$\diagram
          X & \rTo^{g'}   & Z          \cr
   \dTo^{w} &            & \dTo^{f}   \cr
        X_F  & \rTo^{g}  &Y
                                    \enddiagram \quad{\rm with\ } w {\rm\ a\ Whiskering\ given\ by}
           \diagram
         R(F)& \rTo^{i}  &X         \cr
     \dTo^{} &           &\dTo^{w}  \cr
           F & \rTo^{j}  & X_F        
                                    \enddiagram$$ 
We want to define a filler $h:X_F\to Z$, by extending $g'$ along every tree $T$ in $F$.
This is possible since the square is commutative and $f$ is Surjecting.  QED

\bigskip\noindent
{\bf Proposition 3.} (Whiskering, Surjecting) is a weak factorization system.

\medskip\noindent
{\bf Proof:} By the preceding proposition we have 
${\rm Whiskering}\subseteq {}^\dagger{\rm Surjecting}$ and
${\rm Surjecting}\subseteq {\rm Whiskering}^\dagger$.

If $f:X\to Y$ is not in ${\rm Surjecting}$ then there exists some 
$x\in X_0$ and some $a\in Y(f(x),*)$ which is not in the image of
$X(x,*)\to Y(f(x),*)$.  Consider the Whiskering ${\bf s}:{\bf N}\to {\bf A}$ and
the commutative diagram 
$$\diagram
          {\bf N} & \rTo^{x} & X        \cr
   \dTo^{\bf s} &            &\dTo^{f}  \cr
          {\bf A}  & \rTo^a  &Y
                                    \enddiagram$$
for which there is no filler. This shows that
$f\notin{\rm Surjecting}$ implies $f\notin{\rm Whiskering}^\dagger$.
It follows that ${\rm Whiskering}^\dagger={\rm Surjecting}$.
We will show that 
$f\notin{\rm Whiskering}$ implies $f\notin{}^\dagger{\rm Surjecting}$, 
by factoring.
Suppose that $f:X\to Y$ with $f\notin{\rm Whiskering}$; 
then $f=p\circ w$ with some Whiskering $w:X\to X_F$ and 
some Surjecting $p:X_F\to Y$.  
Consider the commutative diagram 
$$\diagram
          X & \rTo^{w}       & X_F        \cr
   \dTo^f &                  &  \dTo^p   \cr
          Y  & \rTo^{\rm id} &Y
                                    \enddiagram$$
If this had a filler $h:Y\to X_F$, then
we would have $p\circ h={\rm id}$. 
This would exhibit $f$ as a ``morphism retract'' of $w$ (a retract in the morphism category).
But we show in the next lemma that this would give us the desired contradiction, finishing our proof.

\medskip\noindent
{\bf Lemma.} Whiskerings are stable  with respect to morphism retract.

\medskip\noindent
{\bf Proof:} First we show that any retract of a rooted tree is a rooted tree.
If $T$ is a rooted tree with root $x_0$ 
and we have the following commutative diagram
with $r\circ s={\rm id}_{T'}$
$$\diagram
    {\bf N} & \rTo^{\rm id}  &{\bf N}    &\rTo^{\rm id} &{\bf N}   \cr
     \dTo^x &                &\dTo^{x_0} &              &  \dTo^x  \cr
          T' & \rTo^{s}       & T         & \rTo^{r}     & T'
                                                    \enddiagram$$
then $T'$ is a rooted tree with root $x$.  This is clear since, for any node $x'$ in $T'$, 
the unique path $\alpha$ from $x_0$ to $s(x')$ gives $r\circ \alpha$ a path
from $x$ to $r(s(x'))=x'$, and there can be no other path in $T'$ from $x$ to $x'$.
 
More generally, consider any commutative diagram
$$\diagram
          X & \rTo^{s'}       & Z     & \rTo^{r'}   & X      \cr
     \dTo^f &                &\dTo^w &            &  \dTo^f  \cr
          Y & \rTo^{s}      & Z_F    & \rTo^{r}  & Y
                                                    \enddiagram$$
such that $w$ is a Whiskering, $r\circ s={\rm id}_Y$ and 
$r'\circ s'={\rm id}_X$.  
The fact that $r\circ s={\rm id}_Y$ implies that $s$ is
injective on nodes and arcs. Also $w$ is injective on nodes and arcs
since it is a Whiskering. This implies that $f$ is injective on
nodes and arcs, since the first square is commutative.

We will describe a rooted forest $F'$ with roots $R'$ a discrete subgraph of $X$, 
such that $Y=X_{F'}$.  

The Whiskering $w$ is given by a rooted forest $F$ whose
$R$ form a discrete subgraph of $Z$.  Let $R'=R\cap X$ as subgraphs of $Z$. 
Then $R'$ is a discrete subgraph of $X$.  For each $x\in R'$, consider the 
tree $T$ in the forest $F$ with root $x_0=s(x)$.  Then $r(T)$ is a retract of $T$,
so $r(T)$ is a tree with root $x=r(x_0)$.  Let $F'=\sum_{x\in R'} r(T)$.
Then $Y=X_{F'}$, and we are done. QED

\bigskip\noindent
{\bf 3. Two Freyd-Kelly factorization systems in graphs.}

The notion of weak factorization system discussed in the preceding section is
a ``weakened'' version of an older notion in category theory.

\medskip\noindent
{\bf Definition:} For morphisms $\ell$ and $r$ in a category ${\cal S}$, 
we say that $\ell\perp r$ ($\ell$ is {\it orthogonal} to $r$) when
$${\rm if} \diagram
          X & \rTo^{f} & A  \cr
   \dTo^{\ell} &           &  \dTo^{r}   \cr
          Y  & \rTo^g   &B
                                    \enddiagram \quad{\rm commutes,\ then}
           \diagram
          X & \rTo^{f} & A  \cr
   \dTo^{\ell} &  \NE^{h} &  \dTo^{r}   \cr
          Y & \rTo^g   &B
                                    \enddiagram\quad {\rm commutes\ for\ {\bf unique}\ } h.$$                   

\medskip
Thus $\ell\perp r$ means that every equation
$rf=g\ell$ has a unique ``solution'' $h$ with $f=h\ell$ and $g=rh$,
so that the original equation is just $r(h\ell)=(rh)\ell$.  

\medskip
Given a class ${\cal F}$ of morphisms in ${\cal S}$ we define
$${\cal F}^\perp=\{r:{\cal F}\perp r\}\quad{\rm and}\quad
{}^\perp{\cal F}=\{\ell: \ell\perp {\cal F}\}\quad{\rm and}\quad
{}^\perp({\cal F}^\perp)\quad{\rm and}\quad({}^\perp{\cal F})^\perp.$$

\medskip\noindent
{\bf Definition:}  A {\it Freyd-Kelly factorization system} in ${\cal S}$ 
is given by two classes
${\cal L}$ and ${\cal R}$ with ${\cal L}^\perp={\cal R}$ and
${\cal L}={}^\perp{\cal R}$, such that for any morphism $c$ in ${\cal S}$,
there exist $\ell\in{\cal L}$ and $r\in{\cal R}$ with $c=r\circ \ell$.

This notion of factorization system
is given in section 2 of Freyd and Kelly [1972].
A Freyd-Kelly factorization system is often just called a {\it factorization system}.
We may use notation $(({\cal L}\ ,\ {\cal R}))$ to indicate a Freyd-Kelly factorization system.

A basic example is the epimorphic, monomorphic
factorization system $(({\cal E}\ ,\ {\cal M}))$ in the category of sets: 
$\{2\to 1\}^\perp={\cal M}$ is the class of injective functions,
and ${}^\perp\{1\to 2\}={\cal E}$ is the class of surjective functions. 
We have ${\cal E}^\perp={\cal M}$ and
${\cal E}={}^\perp{\cal M}$, and every function 
factors as $m\circ e$ with $m\in{\cal M}$ and $e\in{\cal E}$.

\bigskip 

We will describe two interesting Freyd-Kelly factorization systems in Gph
(although we will not need them in this paper).
Our first system is inspired by Stallings [1983].
We will sketch a complete proof, since it involves a nice use of
a ``small object argument'' (for this idea, see Section 2.1 in Hovey [1999], 
for instance).

\medskip
Any pair $\{a',a''\}$ of distinct arcs with a common source node $x$ in a graph $X$
gives a graph morphism $X\to \underline X$,
so that $\{a',a''\}$ becomes a single arc $ a$ in $\underline X$,
and  $\{t(a'),t(a'')\}$ becomes a single node $y$ with $t(a)=y$.
Stallings [1983] called this type of graph morphism  a ``folding''.

\medskip\noindent
Note that the word ``folding'' is used for other types of morphisms in literature on graphs;
see the book Hell and Nesetril [2004], for instance.
We express the Stallings notion of folding as follows.

\medskip\noindent
{\bf Definition:} Let ${\bf V}$ be the graph with three nodes, $0$ and $1'$ and $1''$, and two arcs,
$a'$ and $a''$, with $a'$ from $0$ to $1'$ and  $a''$ from $0$ to $1''$.
Let ${\bf f}:{\bf V}\to {\bf A}$ be the graph morphism taking $a'$ and $a''$ to $a$.
The {\it elementary folding} associated to a graph morphism 
$h:{\bf V}\to X$ is the pushout $X\to \underline X$
of the graph morphisms $h$ and ${\bf f}$.
A {\it Folding} is a (possibly transfinite) composition of elementary foldings.

\medskip\noindent
We include some discussion of transfinite composition of graph morphisms
within the proof of the next proposition.  

The Foldings will form the left class of our first factorization system.
The right class is simpler to define.

\medskip\noindent
{\bf Definition:} A graph morphism $f:X\to Y$ is Injecting if $f:X(x,*)\to Y(f(x),*)$ is a injective
function for all $x\in X_0$.

\medskip\noindent
{\bf Proposition 4.}  Every graph morphism factors as ${\rm Injecting}\circ {\rm Folding}$.

\medskip\noindent
{\bf Proof:} Let $f:X\to Y$ be an arbitrary graph morphism.
We will factor $f$ through the composition of a number of steps.
We essentially show that every object in Gph is small, 
and use a ``small object argument'' (inspired by Section 2.1 in Hovey [1999]).
In fact, if $X$ is infinite we may need a transfinite sequence of steps, 
so we will index our steps by a well-ordered set, an {\it ordinal}.
We view each ordinal as the set of all smaller ordinals 
(see Chapter II, Section 3 in Cohen [1966], for instance).
Every ordinal $\alpha$ has a {\it successor}, defined as $\alpha+1=\alpha\cup\{\alpha\}$.
Let $\Lambda$ be an ordinal so large that there is no 
injective function $\Lambda\to X_1\times X_1$.
We also assume that $\Lambda$ is not the successor of any ordinal, so that
$\lambda\in\Lambda$ implies $\lambda+1\in\Lambda$. 
We view $\Lambda$ as a category
with an object for each element of $\Lambda$ and one morphism from $\lambda$ to $\lambda'$ when
$\lambda\leq\lambda'$.  

Given $f:X\to Y$, we will define a functor $X^\bullet:\Lambda\to {\rm Gph}$ equipped
with natural transformations $f^\bullet$ and $g^\bullet$. 
More precisely, we define $X^\lambda$ and compatible
graph morphisms $f^\lambda:X\to X^\lambda$ and $g^\lambda:X^\lambda\to Y$
by transfinite induction,
assuming that they are defined for all smaller ordinals.
Here each $f^\lambda$ will be a Folding with $f=g^\lambda\circ f^\lambda$.

\smallskip
For the minimal element $0\in \Lambda$, we define $X^0=X$. 
\smallskip
If $g^\lambda$ is not Injecting then we define an elementary folding 
$X^\lambda\to X^{\lambda+1}$
by choosing a node $x$ in $X^\lambda$ having distinct arcs $a',a''\in X^\lambda(x,*)$ such that 
$g^\lambda(a')=g^\lambda(a'')\in Y(g^\lambda(x),*)$.
\smallskip
If $g^\lambda$ is Injecting then we define $X^{\lambda +1}=X^\lambda$ and
$f^{\lambda +1}=f^\lambda$ and $g^{\lambda +1}=g^\lambda$.

\smallskip
For $\lambda$ a limit ordinal (not the successor of any ordinal) 
we define $X^\lambda={\rm colim}_{\alpha<\lambda}X^\alpha$.  The graph morphism
$f^\lambda$, the colimit of Foldings, is called a {\it transfinite composition},
and is a Folding, by our definition.

\smallskip\noindent
Note that if $g^\lambda$ is Injecting, then we will have $X^\lambda=X^{\lambda'}$ 
for all $\lambda'>\lambda$,
and we may say that the {\it $\Lambda$-sequence stabilizes at $\lambda$}.
Let us verify that our $\Lambda$-sequence stabilizes at some $\lambda\in \Lambda$,
so that $f=g^\lambda\circ f^\lambda$ gives our desired factorization.

Each Folding $f^\lambda:X\to X^\lambda$ determines an equivalence relation
$E^\lambda\subseteq X_1\times X_1$ on the arcs of $X$.
So long as $g^\lambda$ is not Injecting,
we have $E^\lambda\subset E^{\lambda+1}$, a strict
inclusion.  This shows that the $\Lambda$-sequence 
constructed above eventually stabilizes, since otherwise we could choose a 
$\Lambda$-parametrized family of elements $p^{\lambda}\in X_1\times X_1$
with $p^{\lambda+1}\in E^{\lambda+1}-E^\lambda$.  This would give an injective function
$\Lambda\to X_1\times X_1$, which is impossible by our assumption about the size of $\Lambda$.
QED

\medskip\noindent
{\bf Proposition 5.} ((Folding, Injecting)) is a Freyd-Kelly factorization system in Gph.

\medskip\noindent
{\bf Proof:}  We have $\{ {\bf f} \}^\perp={\rm Injecting}$, from the definitions.
Since ${\bf f}\in {\rm Folding}$, we have ${\rm Folding}^\perp \subseteq {\rm Injecting}$.
For any class of graph morphisms ${\cal F}$, the class ${}^\perp{\cal F}$ 
is closed under pushouts and under transfinite composition. 
Since ${\bf f}$ generates the Foldings under pushouts
and transfinite compositions, we have ${\rm Folding}\subseteq {}^\perp{\rm Injecting}$.

Finally, we use proof by contradiction to show
${\rm Injecting}^\perp \subseteq {\rm Folding}$.  Suppose
that $g\notin{\rm Folding}$ and $g\in {}^\perp{\rm Injecting}$.
Factor $g=i\circ f$ with $f\in{\rm Folding}$ and $i\in{\rm Injecting}$.
So $f$ is a graph epimorphism.  But the commutative diagram $i\circ f={\rm id}_Y\circ g$
has a (unique) filler $h$; then $i\circ h={\rm id}_Y$ shows that $h$ is a graph monomorphism,
and $h\circ g=f$ shows that $h$ is a graph epimorphism.  So $h$ is a graph isomorphism.
But $f=h\circ g$.  So $g$ is isomorphic to a Folding, which is a contradiction.
QED

\bigskip
We also have a second factorization system, which seems to have interesting connections
to algebraic graph theory.

\medskip\noindent
{\bf Definition:} Let ${\rm Whisker}/{\rm Fold}$ denote the class
of all graphs morphisms which can be factored as a Whiskering followed by a Folding.
Although we will not need it here, we note that
${\rm Whisker}/{\rm Fold}$ can also be described as the graph morphisms
which can be factored as a Folding followed by a Whiskering.

According to the historical sketch given in Boldi and Vigna [2002], the
following basic concept has independently arisen many times in graph theory.  

\medskip\noindent
{\bf Definition:} A graph morphism $f:X\to Y$ is a {\it Covering} when 
$f:X(x,*)\to Y(f(x),*)$ is a bijective function for all $x\in X_0$.

Other names for this concept (and  variants of it) includes {\it divisor}, 
{\it fibration}, {\it equitable partition}, etc; see Boldi and Vigna [2002].

Coverings would seem to play a fundamental role in algebraic graph theory
because of the following.

\medskip\noindent
{\bf Fact:} If $X$ and $Y$ are finite graphs and $f:X\to Y$ is a Covering 
which is surjective on nodes, then the characteristic polynomial of $Y$ 
divides the characteristic polynomial of $X$.

Proofs can be found in Chapter 4 of Cvetkovi\'c, Doob, and Sachs [1978],
and in Section 9.3 of Godsil and Royle [2001].

The connection with algebraic graph theory seems to make the following 
factorization system of graphs especially interesting.
This factorization system was the first that we investigated,
following some remarks by Steve Schanuel.
We state the result without proof here, since we will not need it in this paper:

\medskip\noindent
{\bf Proposition 6.}  $(({\rm Whisker}/{\rm Fold}\ ,\ {\rm Covering}))$ 
is a Freyd-Kelly factorization system in Gph.

\medskip
There is a general principle underlying this.
Let ${\cal S}$ be a suitable category, 
such as the category of Sets or the category Gph;
in particular, ${\cal S}$ is to have colimits.
Let $i_1:Y\to Y+Y$ and $i_2:Y\to Y+Y$ be the two canonical morphisms into the coproduct.
For any morphism $f:X\to Y$, let $Y+Y\to Y+_f Y$ denote the coequalizer of the
two morphisms $i_1\circ f$ and $i_2\circ f$.  
Let $\bar f:Y+_f Y\to Y$ denote the induced morphism.

Let $S$ be a set of morphisms in ${\cal S}$.
Let $\overline S=\{f,\bar f:f\in S\}$.
Let $\overline{\overline S}$ denote the closure of $S$
with respect to the bar operation, all pushouts,  and all transfinite compositions.
Fajstrup and Rosick\'y [2007] show that in the setting of locally presentable categories,
one can guarantee having a factorization system
$$((\ \overline{\overline S}\ ,\ {\overline S}^\dagger\ )).$$

In our example here,
let $S=\{{\bf s}\}$. Then ${\bf f}$ can be identified with $\bar{\bf s}$,
so that $\overline S=\{{\bf s},{\bf f}\}$ and
${\overline S}^\dagger$ is Coverings.
Then $\overline{\overline S}$ is all compositions of Whiskerings and Foldings,
since it is generated by all transfinite compositions of pushouts of ${\bf s}$ and ${\bf f}$ 
(transfinite compositions of
pushouts of ${\bf f}$ are Foldings, and transfinite compositions of
pushouts of ${\bf s}$ are Whiskerings).

\bigskip\noindent
{\bf 4. A model structure on the category of graphs.}

Suppose that  ${\cal S}$ is a category with finite limits and finite colimits.
We take the following definition from Section 7 of Joyal and Tierney [2006].

\medskip\noindent
{\bf Definition:} 
A {\it model structure} on ${\cal S}$ 
is a triple $({\cal C},{\cal W},{\cal F})$
of classes of morphisms in ${\cal S}$ that satisfies

1) ``three for two'': if two of the three morphisms  $a, b, a\circ b$ belong to ${\cal W}$
then so does the third,

2) the pair $(\underline{\cal C} , {\cal F})$ is a weak factorization system 
(where $\underline{\cal C}={\cal C}\cap{\cal W}$),

3) the pair $({\cal C} , \underline{\cal F})$ is a weak factorization system 
(where $\underline{\cal F}={\cal W}\cap{\cal F}$).

\medskip
The morphisms in  ${\cal W}$ are called {\it weak equivalences}.
The morphisms in  ${\cal C}$ are called {\it cofibrations}; and
the morphisms in  $\underline{\cal C}$ are called {\it acyclic cofibrations}.
The morphisms in  ${\cal F}$ are called {\it fibrations}, and
the morphisms in  $\underline{\cal F}$ are called {\it acyclic fibrations}.

Note that, according to Hovey [1999] (page 28),
 ``It tends to be quite difficult to prove that a category admits a model structure.
The axioms are always hard to check.''

\bigskip
Recall from Section 2 that the path graph ${\bf P}_n$ has
nodes $\{0,\ldots,n\}$. 
For $n\geq 0$, the {\it cycle graph} ${\bf C}_n$ is the graph
produced by identifying the nodes $0$ and $n$ of ${\bf P}_n$. 
We have ${\bf C}_0={\bf P}_0={\bf N}$, the graph with one
node and no arcs; and ${\bf C}_1$ is the graph with one node, and one arc with source equal to target.
Let $C_n(X)$ denote the set of graph morphisms from
${\bf C}_n$ to $X$; we may call this the set of
$n$-cycles in $X$.  

\medskip\noindent
{\bf Definition:} A graph morphism $f:X\to Y$ is {\it Acyclic} when
$C_n(f):C_n(X)\to C_n(Y)$ 
is bijective for all $n>0$.

Here we
exclude $n=0$, since we don't want to require that
$f_0:X_0\to Y_0$ is a bijection. 

The Acyclics contain the Whiskerings, and many other 
useful graph morphisms.

\medskip\noindent
{\bf Morphism classes ${\cal W}$ and  ${\cal C}$ and  ${\cal F}$ for 
our Quillen model structure on Gph:} 

Let  ${\cal W}$ be the Acyclics.
Since Whiskerings are Acyclic, we may let $\underline{\cal C}$ be the Whiskerings.
It follows that ${\cal F}$ must be the Surjectings.
Then $\underline{\cal F}$ must be the Acyclic Surjectings.
Finally, we define
${\cal C}$ to be ${}^\dagger{\underline{\cal F}}\ $.
In the next few propositions we show directly that this
does indeed define a model structure on Gph.

\medskip\noindent
{\bf Proposition 7.} The Acyclics satisfy the ``three for two'' property.

\medskip\noindent
{\bf Proof:} This is easy, since Acyclics are defined functorially. 
Consider $h=f\circ g$.  Then $C_n(h)=C_n(f)\circ C_n(g)$ for all $n$.
But the class of bijective functions in the category of sets satisfies the
``three for two'' property.  QED.

We have already shown in Section 2 that $(\underline{\cal C} , {\cal F})$ 
is a weak factorization system.
It remains to do the same for $({\cal C} , \underline{\cal F})$.
We use the following facts, which are easy to verify directly 
from the definition ${\cal C}={}^\dagger{\underline{\cal F}}\ $.

0) Every composition of graph morphisms in ${\cal C}$ is in ${\cal C}$.

1) Every Whiskering is in ${\cal C}$.

2) For any set $I$ and $n>0$, if $A_i\to B_i$ is in ${\cal C}$ for all $i\in I$, 
then $\sum_{i\in I}A_i\to \sum_{i\in I}B_i$ is in ${\cal C}$.

3) For any set $I$ and $n>0$, the graph morphism $\emptyset\to I\times C_n $ is in ${\cal C}$.

4) For any set $I$, the graph morphism $I\times C_n\to C_n$ is in ${\cal C}$.

\noindent In 3) and 4) we view set $I$ 
as a discrete graph and we view graph $I\times C_n$ as a summand of copies of $C_n$.

\medskip\noindent
{\bf Proposition 8.} Every graph morphism $g$ factors as $g=f\circ c$ 
with $c\in {\cal C}$ and $f\in \underline{\cal F}\ $.

\medskip\noindent
{\bf Proof:} Given any graph morphism $g:X\to Y$, we factor $g$ in three steps.

First, we let $C$ be the disjoint union of a copy of ${\bf C}_n$ for each element of $C_n(Y)$
which is not the image of $C_n(g):C_n(X)\to C_n(Y)$.
Let $h:C\to Y$ be the graph morphism which sends each summand cycle of $C$ to its image in $Y$. 
Let $X'=X+C$, and let $g':X\to X'$ denote the inclusion $X\to X+C$, 
and let $f':X'\to Y$ denote the graph morphism $X+C\to Y$ determined by $g:X\to Y$ and $h:C\to Y$.
Then $g=f'\circ g'$, and $g'\in{\cal C}$, and  
$C_n(f'):C_n(X')\to C_n(Y)$ is surjective for all $n>0$.

Next, we let $J=\{(c,n):c\in C_n(Y)\}$, with $j:\sum_J I_c\times {\bf C}_n \to \sum_J {\bf C}_n$
where $I_c$ is the preimage of $c$ for the function $C_n(f'):C_n(X')\to C_n(Y)$.
Also let $k:\sum_J I_c\times {\bf C}_n\to X'$ be the graph morphism which 
sends each summand cycle to the corresponding cycle in $X'$,
and let $\ell: \sum_J {\bf C}_n\to Y$ be the graph morphism which sends 
each summand cycle to the corresponding cycle in $Y$.
Let $g'':X'\to X''$ denote the pushout of $j$ along $k$.
Let $f'':X''\to Y$ be the pushout graph morphism induced by 
$\ell$ and $f'$.
Then $g''\in{\cal C}$ and $f'=f''\circ g''$, and 
$C_n(f''):C_n(X'')\to C_n(Y)$ is bijective for all $n>0$, so that $f''\in {\cal W}$.

Finally, we factor $f''=f'''\circ g'''$ with Whiskering $g''':X''\to X'''$ and
Surjecting $f''':X'''\to Y$, as in Section 2.  Then
$g'''\in{\cal C}$ and $f'''\in{\cal W}\cap {\cal F}\ $.

Thus, $g=c\circ f$ with $c=g'''\circ g''\circ g'$ in ${\cal C}$,
and $f=f'''$ in $\underline{\cal F}$. QED

\medskip\noindent
{\bf Proposition 9.} $({\cal C} , \underline{\cal F})$ is a weak factorization system in Gph.

\medskip\noindent
{\bf Proof:} We have ${\cal C}={}^\dagger\underline{\cal F}$, by definition.
This shows also that $\underline{\cal F}\subseteq{\cal C}^\dagger$.
It remains only to show that ${\cal C}^\dagger\subseteq\underline{\cal F}$.  
But this is easy.  
Consider the Whiskering ${\bf s}:{\bf N}\to{\bf A}$ from Section 2, and the
graph morphisms ${\bf i_n}:0\to{\bf C_n}$,
${\bf j_n}:{\bf C_n}+{\bf C_n}\to {\bf C_n}$ as in 3) and 4) above.
These are all in ${\cal C}$, since each can be lifted against any graph morphism 
in $\underline{\cal F}$.  But if $g\notin \underline{\cal F}$ then we can show failure of
lifting for either ${\bf s}$ or some ${\bf i_n}$ or ${\bf j_n}$. QED

\medskip\noindent
{\bf Corollary:} Our morphism classes ${\cal W}$ and  ${\cal C}$ and  ${\cal F}$
provide a Quillen model structure for the category Gph.

\medskip
The above proofs show how the graph morphisms
${\bf s}:{\bf N}\to{\bf A}$ together with ${\bf i_n}:0\to{\bf C_n}$ and 
${\bf j_n}:{\bf C_n}+{\bf C_n}\to {\bf C_n}$, for $n>0$, generate
our class ${\cal C}$ of Cofibrations.  
This situation is a special case of a general principle in pre-sheaf categories;
see Proposition 7.5 in Joyal and Tierney [2006], for instance.

\bigskip\noindent
{\bf 5. Zeta series and almost isospectral graphs.} 

Ihara zeta functions of graphs are usually discussed in a setting of 
 ``unoriented'' or ``symmetric'' graphs; see Kotani and Sunada [2000], for instance.
We need a version suitable for {\it directed} graphs
(in this section we may refer to objects of our category Gph as directed graphs,
for emphasis). 
There is a nice treatment of zeta series of finite directed graphs in
Section 2 of Kotani and Sunada [2000];  
we will follow them here, but with our own terminology.

\medskip\noindent
{\bf Definition:} A {\it finite graph} is one with finitely many nodes and arcs.
The {\it zeta series} of a finite directed graph $X$ is the formal power series
$$Z(u)={\rm exp}(\sum_{m=1}^\infty c_m{u^m\over m}),$$
where $c_m=|C_m(X)|$ for $m>0$.

See the appendix for some motivation for this definition,
including how it relates to an
Euler product expansion in terms of ``primes''.

\medskip\noindent
{\bf Example:} if $X$ is the graph with one node and $n$ arcs, 
then $c_m=n^m$ and 
$$\sum_{m=1}^\infty c_m{u^m\over m}=\sum_{m=1}^\infty {n^mu^m\over m}=-\log(1-nu)
\quad{\rm so\ that}\quad Z(u)={1\over 1-nu}.$$

\medskip\noindent
{\bf Definition:}  Let $X$ be a finite graph.  Let ${\rm R}X_0$ denote the
real vector space with basis the nodes of $X$.
The {\it adjacency operator $A$ for $X$} is the linear transformation
$A:{\rm R}X_0\to{\rm R}X_0$ determined by 
$$A(x)=\sum_{a\in X(x,*)}t(a)$$
for $x\in X_0$. The {\it characteristic polynomial} of $X$ is defined as 
$a(x)=\det(xI-A)$, the characteristic polynomial of the adjacency operator $A$ for $X$.  
If $X$ has $n$ nodes, then $a(x)$ is a monic polynomial of degree $n$, and
the {\it reversed characteristic polynomial} of $X$ is defined to be
$u^n a(u^{-1})$.

If we totally order the nodes of $X$, then the adjacency operator is represented
by the square matrix $A$ with entry $A_{j,i}$ equal to the number of arcs in $X$ 
from the $i^{th}$ node to the $j^{th}$ node.

Note that the reversed characteristic polynomial of a finite graph $X$ 
has constant term $1$, and is thus a unit in the ring of formal power series  
with integer coefficients.

\medskip\noindent
{\bf Proposition 10.} If $X$ is a finite graph with $n$ nodes 
then the zeta series of $X$ satisfies
$$Z(u)={\rm det}(I-uA)^{-1}={1\over u^n a(u^{-1})}.$$

\medskip\noindent
{\bf Proof:} Let $A$ be any endomorphism of an $n$-dimensional real vector space $V$.
We have ${\rm det}(I-uA)=u^n\det(u^{-1}I-A)$, which proves the second equality.  
One can check the first equality by induction on the dimension of $V$,
since both sides are multiplicative for short exact sequences of vector spaces endowed
with endomorphisms.  Or, when $V$ has a basis of eigenvectors for $A$ 
with eigenvalues $\lambda_1,\ldots,\lambda_n$, we can see the first equality 
from $-\log(1-x)=\sum_k{x^k\over k}$ and
$${\rm exp}(\sum_{m=1}^\infty \sum_{i=1}^n\lambda_i^m{u^m\over m})=
\prod_{i=1}^n \exp(-\log(1-\lambda_i u)=
\prod_{i=1}^n {1\over 1-\lambda_i u}={\rm det}(I-uA)^{-1}.$$
QED

\noindent
{\bf Example (continued):} if $X$ is the graph with one node and $n$ arcs
then $a(x)=x-n$ and $u^1a(u^{-1})=1-nu$, which agrees with $Z(u)={1\over 1-nu}$.

\medskip\noindent
{\bf Proposition 11.} If X and Y are finite graphs 
and $f:X\to Y$ is an acyclic morphism then $Z_X=Z_Y$.

\medskip\noindent
{\bf Proof:} This is clear from the definition, since $|C_m(X)|=|C_m(Y)|$
for all $m>0$.  QED

\medskip\noindent
{\bf Definition:} The eigenvalues of the adjacency operator for $X$ 
may be called the {\it spectrum} of $X$
(even though the adjacency operator is not necessarily a diagonalizable operator).
We say that two finite graphs $X$ and $Y$ are {\it isospectral} 
if they have the same characteristic polynomial.  
We say that $X$ and $Y$ are  {\it almost isospectral} if they have the same
reversed characteristic polynomial.

Loosely speaking, $X$ and $Y$ are almost isospectral if and only if
they have the same non-zero eigenvalues.

\medskip\noindent
{\bf Corollary:} If $X$ and $Y$ are finite graphs with $Z_X=Z_Y$ then $X$ and $Y$ are 
almost isospectral.

\medskip\noindent
{\bf Proof:} This follows immediately from the preceding two propositions. QED

\bigskip\noindent
{\bf An appendix on the zeta series and its Euler product expansion.}

Here is a little history, taken from Thomas [1977],
of how our zeta series for finite directed graphs relates to 
the famous zeta functions from number theory.

\medskip\noindent
{\it The zeta function of Euler and Riemann.} Let $p$ range over the prime numbers.  Then
$$\zeta(s)=\sum_n{1\over n^s}=\prod_p(1-{1\over p^s})^{-1}$$

\noindent
{\it Dedekind's zeta function for algebraic number fields.}  
Let $A$ be the ring of integers in an algebraic number field $K$
(so $K$ is a finite extension over the field ${\rm Q}$ of rational numbers).
Let $N(I)=|A/I|$ for any non-zero ideal in $A$. Then
$$\zeta(s)=\sum_I{1\over N(I)^s}=\prod_P(1-N(P)^{-s})^{-1}$$
where $I$ ranges over the non-zero principal ideals in $A$
and $P$ ranges over the prime ideals in $A$.

\medskip\noindent
{\it A zeta function for algebraic function fields.}
Let $A$ is the ring of integers in an
algebraic function field
(so $K$ is a finite extension over the field ${\rm F}_q(x)$
of rational functions with coefficients in the field ${\rm F}_q$ with $q$ elements).
But here $N(I)=|A/I|=q^{\nu(I)}$ for any non-zero ideal in $A$,
where $\nu(I)$ is the dimension of $A/I$ 
as a finite dimensional vector space over $F_q$, so that
$$\zeta(s)=Z(u)|_{u=q^{-s}} \quad {\rm for} \quad 
Z(u)=\prod_P (1-u^{\nu(P)})^{-1}.$$
The zeta function for a projective variety over a finite field from Weil [1949] is a completed version of this.

\medskip
The form of zeta series that we use in section 5 seems ultimately based
on the following observation.
For a square matrix of a fixed size,
the knowledge of the trace of $A^n$ for all $n$ is equivalent
to the knowledge of the characteristic polynomial of $A$, in that
$${\rm det}(1-uA)^{-1}={\rm exp}(\sum_{n=1}^\infty{\rm Trace}(A^n){u^n\over n}).$$
The proof is just like that of the comparable proposition in section 5.

If $A$ is a permutation matrix 
then the trace of $A^n$ counts the number of fixed points of $A^n$.
This is one of the ideas behind  
the Lefschetz fixed point theorem, the Weil zeta function 
used in the Weil conjectures, and other dynamical zeta functions
such as the Selberg zeta function in Riemannian geometry 
and the Ihara zeta function (see Ruelle [2002] for instance). 

\medskip
In our setting of graphs, we merely use that 
if $A$ is the adjacency matrix of a finite  graph $X$, 
then the trace of $A^n$ counts the number of cycles of length $n$ in $X$. 
This leads to the following Euler product expansion, analogous to the one for algebraic function fields.

The cycle graph ${\bf C}_n$ has nodes $i$ for $0\leq i<n$.
If $m$ divides $n$ then we have a graph morphism $\pi:{\bf C}_n\to {\bf C}_m$
given by sending node $i$ to node $i\bmod m$. 

\medskip\noindent
{\bf Definition:} A cycle $c:{\bf  C}_{km}\to X$ 
is a {\it $k$-multiple} if $c=c'\circ \pi$
for some cycle $c':{\bf C}_m\to X$.
A {\it prime cycle of length $n$ in $X$} is a cycle $c:{\bf C}_n\to X$ 
which is not a $k$-multiple for any $k>1$.
Let us say that two cycles $c,c':{\bf  C}_n\to X$ are {\it shift equivalent} 
if $c'=c\circ \tau^i$ for some $i$,
where  $\tau^i:{\bf  C}_n\to {\bf  C}_n$ is the {\it shift morphism} 
sending node $j$ to node $j+i\bmod n$.  
Let us say that a  {\it prime $P$ in $X$} is an equivalence class of prime cycles in $X$,
and that $\nu(P)$ is the length of the prime $P$.

This makes sense, since shift equivalence is an equivalence relation on $C_n(X)$,
and on prime cycles of length $n$ in $X$.

\medskip\noindent
{\bf Proposition 12.}  The Euler product expansion for the zeta function of a finite graph
is given by
$$Z(u)=\prod_{P}(1-u^{\nu(P)})^{-1}$$ 
where $P$ ranges over all 
primes in $X$ and $\nu(P)$ is the length of $P$.

\medskip\noindent
{\bf Proof:} Let  $\bar c_k$ be the number of primes $P$ of length $k$. 
Then $c_m = \sum_{\{k:k|m\}} k\bar c_k$ and we have
$$\log \prod_{P}(1-u^{\nu(P)})^{-1}=\sum_P\sum_{k=1}^\infty{u^{k|p|}\over k}=
\sum_{k=1}^\infty\sum_{\ell=1}^\infty\sum_{\nu(P)=\ell}{u^{k\ell}\over k}=
\sum_{k=1}^\infty\sum_{\ell=1}^\infty \bar c_\ell{u^{k\ell}\over k}=
\sum_{m=1}^\infty{1\over m}\sum_{\ell|m}\ell \bar{c}_\ell u^m=
\sum_{m=1}^\infty c_m {u^m\over m},$$
which is $\log Z(u)$.  Here equality one and equality four follow from
$$-\log(1-u^v)=\sum_k{u^{kv}\over k}\quad\quad{\rm and}\quad\quad
m=k\ell \iff {\ell\over m}={1\over k}.$$
QED

\medskip\noindent
{\bf Example (continued):} for $X$ the graph with one node and two arcs, we must have
$${1\over 1-2u}=(1-u)^{-2}(1-u^2)^{-1}(1-u^3)^{-2}\cdots=\prod(1-u^k)^{\bar c_k}.$$
This is related to the cyclotomic (necklace) identity; see Dress and Siebeneicher [1989], for instance.

\bigskip\noindent
{\bf An appendix on the Galois connection of a binary relation:}

Suppose we are given a binary relation $\bowtie\ \subseteq{\cal L}\times{\cal R}$.
We write $\ell\bowtie r$ to mean $(\ell,r)\in\ \bowtie$.  Birkhoff [1940]
and Ore [1944] developed the {\it Galois connection} formalism that follows.

\noindent
For $L\subseteq{\cal L}$ and $R\subseteq{\cal R}$ we define
$$L^*=\{r\in{\cal R}:\ \forall \ell\in L,\ \ell\bowtie r\}\quad{\rm and,\ dually,}
\quad {}^*R=\{\ell\in{\cal L}:\ \forall r\in R,\ \ell\bowtie r\}.$$
We define ${}^\bullet L={}^*(L^*)$ 
and say that $L$ is {\it left-closed} iff $L={}^\bullet L$, and dually on the right.

\noindent The following observations are easy to verify:

1) $L_1\subseteq L_2$ implies $L_2^* \subseteq L_1^*$, and dually on the right.

2) $L\subseteq {}^\bullet L$, and dually on the right.

3) $L_1\subseteq L_2$ implies ${}^\bullet  L_1\subseteq {}^\bullet  L_2$, 
and dually on the right.

4) ${}^\bullet {}^\bullet L={}^\bullet  L$, and dually on the right.

5) $L^*=({}^\bullet L)^*$, and dually on the right.

6)  ${}^*R$ is left-closed for any $R$, and dually on the right.

7) ${}^*R_1={}^*R_2$ if and only if $R_1=R_2$, and dually on the right.

8) every left-closed is of the form  ${}^*R$; and dually on the right.

9) the intersection of any collection of left-closeds is left-closed, 
and dually on the right.

\medskip
The above ideas can also be viewed as a basic (contravariant) example of adjoint functors;
see Chapter IV Section 5 in Mac Lane [1971].

\bigskip\noindent
{\bf References:}

\smallskip
\item{[1940]} G. Birkhoff, Lattice Theory, Amer. Math. Soc. Colloquium Publ., 25,
New York, 1940.

\smallskip
\item{[2002]} P. Boldi and S. Vigna, Fibrations of graphs,
Discrete Math., 243 (2002), 21-66.

\smallskip
\item{[2002]} D-C. Cisinski, Th\'eories homotopiques dans les topos,
Journal of Pure and Applied Algebra, 174 (2002), 43-82.

\smallskip
\item{[1966]} P. J. Cohen, Set Theory and the Continuum Hypothesis, W.A. Benjamin, NY, 1966.

\smallskip
\item{[1978]} D. M. Cvetkovi\'c, M. Doob, and H. Sachs, Spectra of Graphs, Academic Press, 1978.

\smallskip
\item{[1983]} A. Dold, Fixed point indices of iterated maps, Invent. math., 74 (1983), 419-435.

\smallskip
\item{[1989]} A. W. M. Dress and C. Siebeneicher, 
The Burnside ring of the infinite cyclic group and 
its relations to the necklace algebra, $\lambda$-rings, and the universal ring of witt vectors.
Adv. in Math., 78 (1989),1-41.

\smallskip
\item{[1995]} W.G. Dwyer and J. Spalinski, Homotopy theories and model categories, 
in Handbook of Algebraic Topology, Editor I. M. James, 73-126, Elsevier, 1995. 

\smallskip
\item{[1999]} E. E. Enochs and I. Herzog,
A homotopy of quiver morphisms with applications to representations,
Canad. J. Math., 51(2) (1999), 294-308.

\smallskip
\item{[2007]} L. Fajstrup and J. Rosick\'y, A convenient category for directed homotopy, arXiv:math/07/08/3937v1

\smallskip
\item{[1972]} P. J. Freyd and G. M. Kelly, Categories of continuous functors, I, 
J. of Pure and Applied Algebra, 2 (1972), 169-191.

\smallskip
\item{[2001]} C. Godsil and G. Royle, Algebraic Graph Theory, Springer-Verlag, New York (2001).

\smallskip
\item{[2004]} P. Hell and J. Ne\v set\v ril, Graphs and Homomorphisms, vol. 28 of
Oxford Lecture Series in Mathematics and its Applications, Oxford University Press, Oxford, 2004.

\smallskip
\item{[1999]} M. Hovey, Model Categories, Amer. Math. Soc., Providence, 1999.

\smallskip
\item{[1991]} A. Joyal and M. Tierney, Strong stacks and classifying spaces, p213-236 in 
Category Theory (Como, 1990), Lecture Notes in Math., 1488, Springer, Berlin, 1991.

\smallskip
\item{[2006]} A. Joyal and M. Tierney, Quasi-categories vs Segal spaces, arXiv:math/06/07/820v1

\smallskip
\item{[2000]} M. Kotani and T. Sunada, Zeta functions of finite graphs,
J. Math. Sci. Univ. Tokyo, 7 (2000), 7-25.

\smallskip
\item{[1989]} F. W. Lawvere, Qualitative distinctions between some toposes of generalized graphs,
In Categories in computer science and logic (Boulder 1987), vol. 92 of Contemp. Math., 261-299, 
Amer. Math. Soc., Providence, 1989.

\smallskip
\item{[1997]} F. W. Lawvere and S. H. Schanuel, Conceptual Mathematics: a first introduction
to categories, Cambridge University Press, Cambridge, 1997.

\smallskip
\item{[1971]} Mac Lane, Categories for the Working mathematician, Graduate Texts in Mathematics,
Springer-Verlag, New York, 1971.
 
\smallskip
\item{[1994]} S. Mac Lane and I. Moerdijk. Sheaves in Geometry and Logic:
a first introduction to topos theory, Universitext, Springer-Verlag, New York (1994)

\smallskip
\item{[1999]} F. Morel and V. Voevodsky, $A^1$-homotopy of schemes, 
Inst. Hautes \'Etudes Sci., Pub. Math. 90 (1999), 45-143.

\smallskip
\item{[1944]} O. Ore, Galois connexions, Trans. Amer. Math. Soc., 55 (1944), 493-513.

\smallskip
\item{[1967]} D. G. Quillen, Homotopical Algebra, Lecture Notes in Mathematics No. 43, Springer-Verlag, Berlin, 1967.

\smallskip
\item{[2002]} D. Ruelle, Dynamical zeta functions and transfer operators, Notices of the Amer. Math. Soc.,
49 no. 8 (2002), 887-895.

\smallskip
\item{[1983]} J. R. Stallings, Topology of finite graphs, Invent. Math., 71 (1983), 551-565. 

\smallskip
\item{[1977]} A. D. Thomas, Zeta-functions: an introduction to algebraic geometry, Pitman, London, 1977.

\smallskip
\item{[1980]} R. W. Thomason, Cat as a closed model category,
Cahiers Topologie G\'eom. Diff\'erentielle, 21 no. 3 (1980), 305-324.

\smallskip
\item{[1949]} A. Weil, Number of solutions of equations in finite fields, Bull. Am. Math Soc., 55 (1949), 497-508.

\bigskip\bigskip\noindent
Terrence P. Bisson, Canisius College, 2001 Main Street, Buffalo, NY 14216 USA

\bigskip\noindent
Aristide Tsemo, College Boreal, 351 Carlaw Avenue, Toronto, Ontario M4K 3M2 Canada

\end